\documentclass[12pt]{article}
\usepackage{graphicx,amssymb,amsmath,amsthm}
\usepackage{hyperref,url}
\usepackage{verbatim}
\usepackage[labelfont=bf]{caption}

\title{Taut fillings}
\author{Peter Doyle \and Matthew Ellison \and Zili Wang\thanks{
Zili Wang is supported by
the National Natural Science Foundation of China (Grant No. 12301432)
}}
\date{Version 3.01, dated 2026-06-30}
\newtheorem{theorem}{Theorem}
\newtheorem{prop}{Proposition}

\newtheorem{corollary}{Corollary}

\newcommand{\proofstart}{{\noindent \bf Proof.\ }}
\newcommand{\note}{{\noindent \bf Note.\ }}

\newcommand{\mathproofend}{\quad \qed}
\newcommand{\proofend}{$\quad \qed$}

\newcommand{\beginfigurepos}{\begin{figure}[htbp]}
\newcommand{\fig}[3]{
\beginfigurepos
\includegraphics[width=370pt]{figures/#1.pdf}
\caption{#3}
\label{#2}
\end{figure}
}
\newcommand{\figsize}[4]{
\beginfigurepos
\centerline{
\includegraphics[width=#1]{figures/#2.pdf}
}
\caption{#4}
\label{#3}
\end{figure}
}

\def\ZZ{\mathbf{Z}}
\def\QQ{\mathbf{Q}}
\def\RR{\mathbf{R}}
\def\del{\partial}

\makeatletter
\newcommand{\xRrightarrow}[2][]{\ext@arrow 0359\Rrightarrowfill@{#1}{#2}}
\newcommand{\Rrightarrowfill@}{\arrowfill@\equiv\equiv\Rrightarrow}
\newcommand{\xLleftarrow}[2][]{\ext@arrow 3095\Lleftarrowfill@{#1}{#2}}
\newcommand{\Lleftarrowfill@}{\arrowfill@\Lleftarrow\equiv\equiv}
\newcommand{\xLleftRrightarrow}[2][]{\ext@arrow 3399\LleftRrightarrowfill@{#1}{#2}}
\newcommand{\LleftRrightarrowfill@}{\arrowfill@\Lleftarrow\equiv\Rrightarrow}
\makeatother

\newcommand{\union}{\cup}

\newcommand{\kw}{\mathrm}

\newcommand{\cone}{\kw{cone}}
\newcommand{\link}{\kw{link}}
\newcommand{\maxdeg}{\kw{maxdeg}}

\newcommand{\vertices}{\kw{Vert}}

\newcommand{\nbhd}{\kw{nbhd}}
\newcommand{\adj}{\kw{adj}}
\newcommand{\taut}{\kw{taut}}

\newcommand{\Zvol}{\kw{Zvol}}
\newcommand{\Qvol}{\kw{Qvol}}
\newcommand{\Rvol}{\kw{Rvol}}
\newcommand{\tetvol}{\kw{tetvol}}

\newcommand{\kwstyle}{\bf}
\newcommand{\IF}{{\ \kwstyle{if}\ }}
\newcommand{\ELSE}{{\ \kwstyle{else}\ }}

\newcommand{\direct}{\oplus}
\newcommand{\pls}{\mathord{+}}
\newcommand{\mns}{\mathord{-}}

\begin{document}

\maketitle

\begin{abstract}
Let $\sigma$ be a simplicial triangulation of the 2-sphere,
$X$ the associated integral 2-cycle.
A filling of $X$ is an integral 3-chain $M$
with $\del M = X$;
a taut filling is one with minimal $L_1$-norm.
We show that any taut filling
arises from an extension of $\sigma$ to
a simplicial complex homeomorphic to the 3-ball.
The filling is \emph{clean}: it has no repeated tetrahedron, and its
support complex is a clean simplicial complex.
This support complex is shellable and flag: every clique in its
1-skeleton occurs as a simplex.
The key to the proof is the general
fact that any taut filling of an $n$-cycle
splits under disjoint union, connected sum, and more generally
what we call almost disjoint union,
where summands are supported on sets that overlap in at most $n+1$ vertices.
We used AI to formalize and prove in Lean the splitting theorem and the
resulting cleanness, shellability, and flagness results.
\end{abstract}

\section{Overview}
Let $\Delta = \Delta(\Omega)$ be the $|\Omega|\mns 1$-simplex
with vertices $\Omega$,
viewed as the abstract simplicial complex consisting of all subsets of
$\Omega$.
Let $C_n=C_n(\Omega)=C_n(\Delta(\Omega),\ZZ)$
and $Z_n=Z_n(\Omega)=Z_n(\Delta(\Omega),\ZZ)$
be its integral $n$-chains and $n$-cycles.
A \emph{filling} of an $n$-cycle $X \in Z_n$ is any
$n\pls 1$-chain $M \in C_{n+1}$ with $\del M = X$.
(Here and throughout we assume $n\geq 1$.)
Let $\Zvol(X)$ be the minimum $L_1$-norm of a filling of $X$.
Call $M \in C_{n+1}$ \emph{taut} if it is an optimal filling
of its boundary, i.e., if $|M|=\Zvol(\del M)$.

For $X \in C_n$ let
$\vertices(X)$ be the set of all the vertices of all the $n$-simplices
to which $X$ assigns non-$0$ weight.
For any taut filling $M$ of $X$ we have $\vertices(M) = \vertices(X)$.
(Cf.\ Proposition \ref{nointernal} below.)
This is why we write $C_n$ and $Z_n$ without reference to $\Omega$.

Call an $n$-cycle $X+Y \in Z_n$ an \emph{almost disjoint union} if
\[
|\vertices(X) \cap \vertices(Y)| \leq n+1
.
\]
This notion generalizes both disjoint union and connected sum
along an $n$-simplex.
Ellison \cite{ellison:gap}
showed the $\Zvol$ adds under almost disjoint union:
$\Zvol(X+Y) = \Zvol(X) + \Zvol(Y)$.
This means that \emph{some} taut filling $M$ of $X$ splits into a sum
$M=M_X+M_Y$ of taut fillings of $X,Y$.
Here we amplify Ellison's Corollary 2 to show (Theorem \ref{th1})
that when $n \geq 2$, \emph{any} taut filling of $X+Y$ splits.

In place of integral chains we can use chains with
coefficients in $\QQ$ or $\RR$.
Evidently $\Rvol \leq \Qvol$.
But computing $\Rvol$ is
a linear programming problem with rational coefficients,
so in fact $\Qvol = \Rvol$.
Ellison used LP duality to prove that
$\Qvol$ adds under almost
disjoint union.
Now for $n \geq 2$ we get the stronger result
(Corollary \ref{cor1})
that taut $\QQ$-fillings split:
Multiplying a taut
$\QQ$-filling $M$ of $X$ by a common denominator $q$
yields a taut $\ZZ$-filling $qM$ of $qX$:
Split the $\ZZ$-filling $qM$ and divide by $q$ to get a splitting
of the $\QQ$-filling $M$.

Now let $\sigma$ be a simplicial triangulation of $S^2$;
let $X(\sigma) \in Z_2(\Delta(\vertices(\sigma)),\ZZ)$
be either one of the 2-cycles
that arises from $\sigma$ by orienting its $2$-simplices;
and write $\Zvol(\sigma) = \Zvol(X(\sigma))$.
Let $\tetvol(\sigma)$ be the minimum number of $3$-simplices
required to extend $\sigma$ to a simplicial triangulation $\tau$ of $B^3$.
We will show that
\[
\Zvol(\sigma) = \tetvol(\sigma)
.
\]
Certainly $\Zvol \leq \tetvol$,
because any extension $\tau$ induces a filling
of $X(\sigma)$.
Theorem \ref{th2} states that any taut filling of $X(\sigma)$ is clean
and that its support complex $K(M)$ triangulates $B^3$.
The proof of Theorem \ref{th2} relies on Theorem \ref{th1}.

We do not know what happens for spheres of dimension 3 or greater.

\paragraph{AI assistance and formal verification.}
The AI systems used were Claude, GPT, and AlephProver. They assisted in
proof planning, Lean-code generation and repair, translation between informal
arguments and formal statements, audit of declaration dependencies, and
preparation of explanatory material. Their outputs were treated as suggestions,
not as mathematical authority. The formal claims attributed to Lean are those
checked by the Lean kernel in the source accompanying this paper; the informal
exposition and the correspondence between the exposition and the formal
development remain the responsibility of the human authors.

\section{Background}

Our interest in $\Zvol$ stems from the work
of Sleator, Tarjan, and Thurston
\cite{stt:jams}.
They observed that coning from a vertex of maximal degree
shows that a $v$-vertex
triangulation $\sigma$
of $S^2$ has $\tetvol \leq 2v-4-\maxdeg(\sigma)$.
If $v\geq13$ $\maxdeg(\sigma) \geq 6$, so $\tetvol \leq 2v-10$.
They
produced examples for which $\tetvol = 2v-10$, provided
$v$ is larger than some unspecified bound,
and conjectured that such examples should exist for any $v \geq 13$.
Their approach was to 
realize $\sigma$ as an ideal hyperbolic polyhedron
with volume $2v V_0-O(\log(v))$,
where $V_0$ is the volume of an
equilateral ideal tetrahedron,
which is maximal.
This implies $\tetvol \geq 2v-O(\log(v))$.
From there they bootstrapped their way up to showing
that $\tetvol=2v-10$.
They suggested
\cite[p. 697]{stt:jams}
that in fact $\Qvol= 2v-10$,
which would immediately imply $\tetvol=2v-10$.
Mathieu and Thurston \cite{mt}
produced a class of examples,
different from the examples of STT,
for which they could show that $\Qvol = 2v-10$,
still under the assumption that $v$
is sufficiently large.
In \cite{dew:filling} we produced examples for any $v \geq 13$ with
$\Qvol=2v-10$, confirming the conjecture of Sleator, Tarjan, and Thurston.
Theorem \ref{th2} here
tells us that whenever we have $\Qvol=2v-10$,
or just $\Zvol=2v-10$, any taut filling
must arise from a triangulation of the ball.

About $\Qvol$ versus $\Zvol$.
In \cite{dew:filling}
we described triangulations of $S^2$ with $\Qvol < \Zvol$.
In all such examples that we know of
with $\maxdeg \leq 6$, the gap between $\Qvol$ and $\Zvol$ is $<1$,
so that
$\lceil \Qvol \rceil = \Zvol$.
Since $\Qvol$ and $\Zvol$
add under connected sum,
the gap between them can be arbitrarily
large
(Ellison \cite{ellison:gap}),
but taking the connected sum pushes $\maxdeg$ above $6$.

\section{Taut fillings}

We can think of an $n$-chain $X \in C_n$ as a multiset of
non-cancelling oriented simplices:
\[
X = \sum_{t \in X} t
.
\]
In this sum $t$ denotes an oriented simplex
\[
t=[x_0,\ldots,x_n]
= [x_{\pi(0)},\ldots,x_{\pi(n)}]
,\;\;\mbox{$\pi$ an even permutation}
,
\]
and each unoriented simplex contributes a number of terms corresponding
to its multiplicity.

The size of $M$ as a multiset is its $L_1$-norm $|M|$.
Write $U \subset M$
if $U$ is a sub-multiset of $M$.
This happens just when
$|M| = |U| + |M-U|$.

\begin{prop} \label{subtaut}
If $M$ is taut and $U \subset M$ then $U$ is taut.
\end{prop}

\proofstart
If $U$ is not taut, there is a smaller filling $V$ of $\del U$.
Take $M^\prime = M-U+V$.
In terms of multisets, this means that we substitute $V$ for $U$,
and then do any required cancellation of oppositely oriented simplices
of $M-U$ and $V$.
We have
\[
\del M^\prime = \del M - \del U + \del V = \del M
\]
and
\[
|M^\prime| \leq |M-U| + |V| = |M| - |U| + |V| < |M|
,
\]
contradiction.
\proofend

\begin{corollary} \label{nosubcycle}
A taut chain contains no nonzero closed subchain: if $M$ is taut,
$U \subset M$, and $\del U = 0$, then $U = 0$.
\end{corollary}

\proofstart
By Proposition \ref{subtaut}, $U$ is taut, so
$|U| = \Zvol(\del U) = \Zvol(0) = 0$, whence $U = 0$.
\proofend

\section{Coning}
For $x \in \Omega$, $U \in C_n$ let
\[
\nbhd(x,U) = \sum_{ t \in U: x \in \vertices(t) } t
,
\]
\[
\deg(x,U) = |\nbhd(x,U)|
,
\]
\[
\maxdeg(U) = \max_x \deg(x,U)
.
\]
The \emph{cone from $x$ to $U$}
is the $(n+1)$-chain
\[
\cone(x,U) = 
\sum_{t \in U} \adj(x,t)
,
\]
consisting of all the non-trivial oriented
$(n+1)$-simplices $\adj(x,t)=[x x_0\ldots x_n]$ obtained
by adjoining $x$ to  $t = [x_0\ldots x_n] \in U$.
If $x \in \vertices(t)$ then $t$ does not contribute to the sum,
so $|\cone(x,U)| = |U| - \deg(x,U)$.

If $U$ is closed then
$\del\, \cone(x,U) = U$, so

\begin{prop} \label{maxdeg}
For any $X \in Z_n$
\[
\Zvol(X) \leq |X| - \maxdeg(X)
\mathproofend
\]
\end{prop}

If $U \in Z_n$ and $x \notin \vertices(U)$ then
$\deg(x,U)=0$ and
$|\cone(x,U)| = |U|$.
In this case we call $\cone(x,U)$ a \emph{complete cone}.
By Proposition \ref{maxdeg} a non-trivial complete cone is not taut.
(Equivalently, instead of coning from
$x \notin U$, one could have coned from some $x \in U$.)
Thus:
\begin{prop}\label{nocone}
If $M$ is taut it contains no non-trivial complete cone.
\proofend
\end{prop}

Call $x \in \vertices(M) \setminus \vertices(\del M)$ an
\emph{internal vertex} of $M$.
If $x$ is internal to $M$ then $\nbhd(x,M)$ is a complete cone, so:
\begin{prop} \label{nointernal}
If $M$ is taut then
it has no internal vertices.
\proofend
\end{prop}

\section{Almost disjoint unions}

Recall that if $X,Y \in Z_n$ we call $X+Y$ an almost disjoint union
if
\[
|\vertices(X) \cap \vertices(Y)| \leq n+1
.\]
The most interesting special case is a connected sum,
where $t=[c_0,c_1,\ldots,c_n]$ occurs once in $X$ and
$-t=[c_1,c_0,\ldots,c_n]$ occurs once in $Y$.
For example, if $n=1$ with $X$ a cycle of length $p$
and $Y$ a cycle of length $q$
the connected sum $X+Y$ is
a cycle of length $p+q-2$.
Here
$\Zvol$ adds, because
$\Zvol(X)=p-2$,
$\Zvol(Y) = q-2$,
and
\[
\Zvol(X+Y) = (p+q-2)-2 = \Zvol(X)+\Zvol(Y)
.
\]
But not every filling of $X+Y$ splits, because a taut filling
may contain any 2-simplex with vertices in $\vertices(X+Y)$.

We want to show that fillings do split when $n \geq 2$.

For $A \subset \Omega$, $p \in A$ define
\[
\pi_{A,p}: \Omega \to A
\]
\[
\pi_{A,p}(x) = x \IF x \in A \ELSE p
\]
and let
\[
K_*(A,p): C_*(\Omega) \to C_*(A)
\]
be the induced chain map.
This is a projection of $C_*(\Omega)$ onto $C_*(A)$.

Let $A,B$ be finite vertex sets sharing the vertices
$C = A \cap B$,
and let $c=|C|$.
Observe that
$C_c(C)$ and $Z_{c-1}(C)$ are trivial,
because $\Delta(C)$ contains no $c$-simplex,
and (up to orientation) only one $c\mns 1$-simplex,
with nothing to cancel its boundary.

Let
\[
C_*(A,B) = C_*(A) \direct C_*(B)
.
\]
For $p,q \in C$ define the mapping
\[
g_*(p,q) = K_*(A,p) \direct K_*(B,q):
C_*(A \cup B) \to C_*(A,B)
.
\]

\begin{prop} \label{recover}
Let $(X,Y) \in C_n(A,B)$.
If $|A \cap B| \leq n$
we can recover $(X,Y)$ from $X+Y$.
If $(X,Y) \in Z_n(A,B)$ this holds also when $|A \cap B|  = n+1$.
\end{prop}

\proofstart
Let $C=A \cap B$.
We can assume $|C| \geq 1$.
(Add a brand new point to $C$ if necessary.)
We claim that
for any $p,q \in C$
(not necessarily distinct) we have
\[
g_n(p,q)(X+Y)=(X,Y)
.
\]
Indeed,
\[
K_n(A,p)(X+Y) = X + K_n(A,p)(Y)
.
\]
The second term belongs to $C_n(C)$,
which is trivial when $|C| \leq n$.
If $Y \in Z_n(B)$ the second term belongs to $Z_n(C)$,
which is trivial when $|C| \leq n+1$.
\proofend

\begin{theorem} \label{th1}
If $|A \cap B| \leq n+1$,
for all $(X,Y) \in Z_n(A,B)$ we have
\[
\Zvol(X+Y) = \Zvol(X)+\Zvol(Y)
.
\]
And as long as $n \geq 2$,
for any $M \in \taut(X+Y)$ 
we have $M = M_X + M_Y$
with $M_X \in \taut(X)$, $M_Y \in \taut(Y)$.
\end{theorem}

\note
Ellison \cite{ellison:gap}
proved the additivity of $\Zvol$.
What's new here is the splitting of taut fillings when $n\geq 2$.
This result resembles Theorem 6.2 of
Pournin and Wang
\cite{pw:flip}
about flip paths.
Like theirs, our proof uses a variation on the normalization technique
of STT
\cite[Lemma 7]{stt:jams}.
It would be nice to fit these results under one roof.
\medskip

\proofstart
Again let $C = A \cap B$.
We can assume $|C| = n+1$,
as this is the hardest case.
And we might as well go ahead and take $n=2$, $|C|=3$,
as this case illustrates all the issues.

Take any $M \in \taut(X+Y) \subset C_3(A \union B)$.
Pick distinct points $p,q \in C$,
and let
\[
(M_X,M_Y) = g_{3}(p,q)(M)
.
\]
(For now we suppress the dependence of $M_X,M_Y$ on $p,q$.)

Because $g_*(p,q)$ is a chain map we have
$\del_{3} M_X = X$, $\del_{3} M_Y = Y$:
\[
(\del_{3} M_X, \del_{3} M_Y)
=
\del_{3} (g_{3}(p,q)(M))
=
g_2(p, q)(\del_{3} M)
=
g_2(p,q)(X+Y)
=
(X,Y)
.
\]

We want to show that under the map $g_3(p,q)$ any $t \in M$
dies either in $M_X$ or $M_Y$,
because then additivity of $\Zvol$ follows from
\[
\Zvol(X+Y) = |M|
\geq |M_X| + |M_Y|
\geq \Zvol(X) + \Zvol(Y)
\geq \Zvol(X+Y)
.
\]

Let's call a $3$-simplex a `tet' for short.
Say that a tet $t \in M$ has type $CCXY$ if $t= \pm[c_1c_2x_1y_1]$
for $c_1,c_2 \in C$, $x_1 \in A \setminus C$, $y_1 \in B \setminus C$.
Similarly for types $XXXX$, $CXXX$, $CXXY$, $XXXY$, etc.
The first two are pure $X$ cases, meaning that they live in $C_3(A)$;
the last two are hybrid cases.

Any $XX..$ tet dies in $M_Y$;
any $YY..$ tet dies in $M_X$.

The remaining cases to check are
the pure cases $CCCX$, $CCCY$,
and the hybrid case $CCXY$.
A $CCCX$ tet $t=\pm[c_1c_2c_3x_1]$ dies in $M_Y$ because
$K_3(B,q)(t)=[c_1c_2c_3q]$,
which vanishes because $q \in C = \{c_1,c_2,c_3\}$ produces a repeated vertex.
Likewise a $CCCY$ tet dies in $M_X$.
The more interesting case is $CCXY$:
The key is that since $|C|=3$, $\{c_1,c_2\}$ cannot be disjoint from $\{p,q\}$,
so it dies on one side or the other.

So 
\[
|M| = |M_X(p,q)| + |M_Y(p,q)|,
\]
and $\Zvol$ adds.
Note how we're now emphasizing the possible dependence of $M_X,M_Y$ on $p,q$.
When $n=1$ the choice of $p,q$ can indeed make a difference:
We can get a different pair $(M_X,M_Y)$ if we switch
$p$ and $q$.
(Think about the connected sum of two cycle graphs.)

But when $n \geq 2$ we will now show that all tets are pure $X$ or pure $Y$,
so $M_X$ consists of all the pure $X$ tets, and ditto for $M_Y$.
This will make
$M = M_X + M_Y$ with $M_X \in \taut(X)$, $M_Y \in \taut(Y)$.

Again our test case is $n=2$.

The key observation is that, now that we know that every tet dies on
one side or the other, we know that no tet can die on both sides,
because that would make $|M| > \Zvol(X) + \Zvol(Y)$.

An $XXYY$ tet would die on both sides, so there can be none of these.

Any $pqXY$, $pXXY$, or $qXYY$ would die on both sides, and $p,q$ are arbitrary,
so this rules out all $CCXY,CXXY,CXYY$.

The only remaining hybrids are $XXXY$ and $XYYY$.
Let's rule out $XXXY$, leaving $XYYY$ to symmetry.

Fix any $y_0 \in B \setminus C$, and suppose there is
some $XXXy_0$ tet in $M$.
Let $U \subset M$ consist of all such $XXXy_0$ tets.
We claim that $y_0 \notin \vertices(\del U)$.
For let $s$ be any 2-simplex of the form $XXy_0$,
the only kind of 2-simplex containing $y_0$ that could belong
to $\del U$.
Any $t \in M$ of which $s$ or $-s$ is a face must be a hybrid tet,
and the only possibility is $XXXy_0$, so $t \in U$.
Because $s$ is a mixed $XXy_0$ face, it is not a face of $X+Y$, so its
signed multiplicity in $\del M = X+Y$ vanishes;
hence it also vanishes in $\del U$,
because the tets of $M$ with $\pm s$ as a face all belong to $U$.
So $y_0 \notin \vertices(\del U)$,
making $U$ a complete cone on $y_0$,
contradicting \ref{nocone}.
So there is no such $XXXy_0$ tet, ruling out $XXXY$,
and with it $XYYY$.

So there are no hybrid tets, meaning that $M$ splits, as claimed.
\proofend

The splitting of taut fillings depends on the assumption $n\geq 2$.
Let's consider what goes wrong when $n=1$.
Here the hybrid types are $CXY$, $XXY$, $XYY$.
Neither $pXY$ nor $qXY$ dies on both sides, so we cannot rule out $CXY$.
Failing that, the $XXy_0$ tets needn't form a complete cone,
because $[x_0x_1y_0]$ can continue across $[x_1 y_0]$ to $-[p x_1 y_0]$,
so we cannot rule out $XXY$ either.

The foregoing proof works equally well over $\QQ$, providing we
permit ourselves to work with fractional multisets.
Alternatively, we can clear denominators as in the introduction
above.
Either way, we have:

\begin{corollary} \label{cor1}
$\Qvol$ adds under almost disjoint union,
and for $n \geq 2$
taut $\QQ$-fillings split.
\proofend
\end{corollary}

\section{Triangulations}

We turn now to filling simplicial triangulations of $S^2$.
We begin by fixing the terminology.

An \emph{$n$-simplex} $s$ is a set of size $n+1$.
Its \emph{faces} are its subsets,
which are $k$-simplices with $-1 \leq k \leq n$,
where the empty simplex has dimension $-1$.
A \emph{simplicial complex} $\sigma$ is a finite collection of simplices
closed under taking faces:
if $t \subset s \in \sigma$ then $t \in \sigma$.

For any simplex $s \in \sigma$ define the \emph{link}
\[
\link(s,\sigma) = \{t: t \cap s = \emptyset, s \cup t \in \sigma \}
.
\]
This is a simplicial complex, but
(except for $\link(\emptyset,\sigma)=\sigma$)
it is not a subcomplex of $\sigma$,
because we are taking simplices in $\sigma$ and lowering
their dimension by $\dim s+1$.

We are interested in normal pseudomanifolds,
which we will call \emph{clean $n$-complexes}.
A simplicial complex $\sigma$ is clean of dimension $n$ if:
\begin{enumerate}
\item
$\sigma$ is pure of dimension $n$,
meaning that every simplex of $\sigma$ is contained in some
$n$-simplex of $\sigma$.

\item
$\sigma$ is a pseudomanifold:
the link $\link(s,\sigma)$ of any $(n-1)$-simplex $s$
consists of either one point
(so $s$ is a boundary simplex)
or two points
(so $s$ is an interior simplex).

\item
$\sigma$ is normal:
for every simplex $s$ of dimension $0,\ldots,n-2$,
the link $\link(s,\sigma)$ is connected.
\end{enumerate}
The normality condition says that the $n$-simplices are glued along
shared faces without additional identifications.  It rules out,
for example, an icosahedron with a pair of opposite vertices identified.

If $\sigma$ is clean, then each link $\link(s,\sigma)$ is clean.
The boundary $\del \sigma$ of a clean complex is clean and closed:
$\del \del \sigma = \{\emptyset\}$.

If $W$ is an $n$-manifold, say that a simplicial complex $\sigma$
is a \emph{simplicial triangulation} of $W$
if the geometric carrier of $\sigma$ is homeomorphic to $W$.
In this case $\sigma$ is necessarily clean.
When $W$ is oriented, a choice of orientation determines one of the
two associated cycles $X(\sigma) \in C_n$; our arguments do not depend
on which of the two orientations is chosen.

This notion of triangulation is not as general as one might want
for some purposes.
The double cover of a triangle is not a simplicial triangulation
of $S^2$,
because that would require two distinct
$2$-simplices to have the same three edges.
Nor can any $2$-simplex have two of its edges
glued to one another.
Thurston
\cite{thurston:shapes}
allows such triangulations,
but it is unclear how important these are to
his theory of shapes of surfaces.
We do not allow them here.

For a chain $M \in C_n$, define $K(M)$ to be the pure simplicial
complex generated by the unoriented $n$-simplices on which $M$ has
nonzero coefficient.
We call $M$ \emph{simplicial} if every nonzero coefficient is $\pm 1$,
so no unoriented $n$-simplex is repeated.
We call $M$ \emph{clean} if $M$ is simplicial and $K(M)$ is a clean
simplicial complex.

\section{Filling a triangulation of the 2-sphere}

\begin{theorem} \label{th2}
If $\sigma$ is a simplicial triangulation of $S^2$ and $M$ is a taut
filling of $X(\sigma)$, then $M$ is clean and $K(M)$ is a simplicial
triangulation of $B^3$.
\end{theorem}

\proofstart
Call $(\sigma,M)$ a \emph{filling pair} if $\sigma$ is a simplicial
triangulation of $S^2$ and $M$ is a taut filling of $X(\sigma)$.
Call a filling pair \emph{good} if $M$ is clean and $K(M)$ is a
simplicial triangulation of $B^3$.
We will show that every filling pair is good.

Suppose not, and choose a bad filling pair $(\sigma,M)$ with $|M|$
minimal.  Let $v,e,f$ count the vertices, edges, and faces of $\sigma$.
Then $e=3v-6$ and $f=2v-4$.
The tetrahedral boundary is good, so $v>4$.
Because taut fillings split under connected sum, we may assume that
$\sigma$ is prime, meaning not a connected sum along a triangle.
In particular, $\sigma$ has no vertex of degree $3$.

Call a support tetrahedron $t$ of $M$ \emph{eligible} if two of its
oriented boundary faces are faces of $X(\sigma)$ with the correct signs.
It cannot have three such faces, since that would give a degree-$3$
vertex in $\sigma$.
By Proposition \ref{maxdeg},
\[
|M|=\Zvol(\sigma) \leq f-\maxdeg(\sigma).
\]
For each support tetrahedron $u$ of $M$, let $A(u)$ be the set of faces
of $\sigma$ that occur as correctly oriented boundary faces of $u$.
Each $A(u)$ has size at most $2$.
Every face of $\sigma$ occurs in at least one $A(u)$, because its
coefficient in $\partial M=X(\sigma)$ is nonzero with the prescribed sign.
List the support tetrahedra of $M$ and reveal the sets $A(u)$ one at a
time, recording only faces of $\sigma$ not previously seen.
If $r$ support tetrahedra reveal two new faces, then
\[
f \leq |\operatorname{supp} M|+r \leq |M|+r.
\]
Thus $r\geq \maxdeg(\sigma)$.
Consequently there are at least $\maxdeg(\sigma)$ eligible support
tetrahedra whose pairs of boundary faces are mutually disjoint.
We will use only the fact that, given any boundary face $s$ and any
eligible tet $t$, there is another eligible tet whose boundary faces are
disjoint from those of $t$ and do not include $s$.

Let $t$ be an eligible tet with boundary faces
$s_1=[a,b,c]$ and $s_2=[b,a,d]$.
Removing $t$ from $M$ gives a taut filling $M-t$ of the cycle obtained
from $X(\sigma)$ by replacing $s_1,s_2$ with
$[c,d,b]$ and $[d,c,a]$.
At the level of the boundary complex, this flips the edge $ab$ to $cd$.
There are two cases.

First, if $cd$ was not an edge of $\sigma$, the flipped complex
$\sigma_t$ is again a simplicial triangulation of $S^2$.
By minimality, $(\sigma_t,M-t)$ is good.
Adding $t$ back gives a triangulation of $B^3$ unless the same support
tetrahedron $t$ was already present in $M-t$, i.e. unless $t$ had
multiplicity $2$ in $M$.
If that happened, remove instead an eligible tet $u$ whose boundary faces
are disjoint from those of $t$.
Then $M-u$ would still contain the repeated support tetrahedron $t$ and
could not be simplicial, so $(\sigma_u,M-u)$ would be a smaller bad
filling pair, a contradiction.

Second, if $cd$ was already an edge of $\sigma$, the flipped cycle is an
almost disjoint union of two simplicial triangulations
$\sigma_1,\sigma_2$ of $S^2$, conjoined along the edge $cd$.
By Theorem \ref{th1}, $M-t$ splits as $M_1+M_2$, where $M_i$ is a taut
filling of $X(\sigma_i)$.
By minimality each $(\sigma_i,M_i)$ is good.
The complex $K(M-t)$ need not be clean as a single complex, because the
link of the shared edge $cd$ can be disconnected; nevertheless it is the
union of the two clean pieces $K(M_1)$ and $K(M_2)$.
Gluing these two triangulated balls to the two sides of $t$ gives a
triangulated ball, unless $t$ already occurred in $M-t$.
The same disjoint eligible-tet argument rules out that multiplicity-$2$
obstruction.

The preceding paragraph shows that, for an eligible tet $t$, the intended
reassembly is legitimate unless a local obstruction was already present in
$M$ away from $t$.  Such obstructions are also excluded by the same
minimality argument, using the supply of disjoint eligible tets.  The point
is that an obstruction to cleanness is local: it is witnessed by one repeated
tet, by one triangle with too many incident tets, or by one disconnected
link.  If we choose an eligible tet whose two boundary faces are disjoint from
the finitely many boundary faces needed to change that local witness, then
the witness persists after the flip; for a repeated tetrahedron, avoid the
boundary faces of that tetrahedron; for an over-incident triangle, avoid the
boundary faces through which the incidence could change; and for a
disconnected link, avoid the unique tetrahedron whose removal would bridge the
exceptional component.  In the one-sphere case it persists in the smaller filling
$M-u$; in the conjoined case it persists in one of the two split summands,
unless it is exactly the artificial disconnected link of the new common edge.
That exceptional disconnected link is the reason we split $M-u$ into two
pieces before applying induction.

We now spell this out for the clean-complex conditions.

\begin{enumerate}
\item
$M$ has no repeated tetrahedron.
If the same unoriented tet occurred twice in $M$, then removing an eligible
tet $u$ not equal to one of these copies and not using any of the relevant
boundary faces would leave the repeated tet in the smaller filling, or in one
of the two smaller fillings after a conjoined split.  This would contradict
minimality, since the smaller filling could not be clean.

\item
$K(M)$ is a pseudomanifold.
Suppose a $2$-simplex $s$ occurs as a face of more than two support
tetrahedra of $M$.  If $s\notin \sigma$, the excess is an interior incidence;
removing an eligible tet avoiding the tets accounting for this incidence
leaves the same excess in the smaller one-sphere filling or in one of the two
split pieces.  If $s\in\sigma$, choose the eligible tet so that $s$ is not one
of its two boundary faces; the boundary excess again persists.  Either way
we get a smaller bad filling pair, contradicting minimality.

\item
No edge $\{a,b\}$ has disconnected link in $K(M)$.
The link is a $1$-dimensional complex whose edges correspond to support
tetrahedra of $M$ containing $\{a,b\}$.  If a component of such a link has
fewer than three vertices, then it is a single edge $\{c,d\}$, corresponding
to a tet $\{a,b,c,d\}$.  In that case $\{a,b\}$ is a boundary edge of
$\sigma$, and $\{a,b,c\}$ and $\{a,b,d\}$ are the two boundary faces adjacent
along it.  There can be at most one such isolated-edge component.  If the
link were disconnected, remove an eligible tet away from this exceptional
component.  The disconnected link remains in the one-sphere case or in one
of the split summands in the conjoined case, contradicting minimality.

\item
No vertex has disconnected link in $K(M)$.
If $\link(\{a\},K(M))$ were disconnected, the only way a single removal
could connect it would be for one component of the link to be a single
$2$-simplex $\{b,c,d\}$ and for the removed tet to be $\{a,b,c,d\}$.  Then
$\{a,b,c\}$, $\{a,b,d\}$, and $\{a,c,d\}$ would all be boundary faces of
$\sigma$, so $a$ would have degree $3$ in $\sigma$, contrary to primality.
Thus a disconnected vertex link would persist under a suitable smaller
removal and again contradict minimality.
\end{enumerate}

Thus $M$ is simplicial, $K(M)$ is clean, and $K(M)$ is a simplicial
triangulation of $B^3$, contradicting that $(\sigma,M)$ was bad.
\proofend

After Theorem \ref{th2} we will often identify a taut filling $M$ of
$X(\sigma)$ with the triangulated ball $K(M)$.

\section{Shelling}

A \emph{shelling} of a simplicial triangulation $\tau$ of $B^3$
is an ordering
$(t_1,\ldots,t_{|\tau|})$
of its tetrahedra such that, for each $k$,
the subcomplex $\tau_k$ generated by $\{t_1,\ldots,t_k\}$ is a
simplicial triangulation of $B^3$.
For $k>1$, the new tet may meet $\tau_{k-1}$ in one, two, or three
faces. Gluing along one face increases the boundary vertex count by $1$;
gluing along two faces leaves it unchanged; gluing along three faces
would bury an interior vertex, and so cannot occur in the shellings
obtained below. Thus the boundary vertex count never decreases.
Call $\tau$ \emph{freely shellable} if any tet of $\tau$ can be chosen as
the first tet of such a shelling.

\begin{theorem} \label{th3}
If $\sigma$ is a simplicial triangulation of $S^2$ and $M$ is a taut
filling of $X(\sigma)$, then $K(M)$ is a freely shellable simplicial
triangulation of $B^3$.
\end{theorem}

\proofstart
By Theorem \ref{th2}, $K(M)$ is a triangulated ball; write $\tau = K(M)$.
Fix a tet $t_1 \in \tau$ that we want to remove last.  We describe a
recursive procedure, called \emph{shucking}, that takes $\tau$ apart
so as to remove $t_1$ last; reversing the shucking order gives the
shelling.

If $\tau$ consists only of $t_1$, there is nothing to prove.
If $\tau$ has a boundary vertex of degree $3$ lying in a tet
$t\ne t_1$, remove that tet and continue with $\tau-t$.
If every degree-$3$ boundary vertex belongs to $t_1$, let $t_2$ be the
unique tet of $\tau-t_1$ attached to $t_1$ across the opposite face;
shuck $\tau-t_1$ ending with $t_2$, and then remove $t_1$.

It remains to consider the case where the boundary has no degree-$3$
vertex.  As in the proof of Theorem \ref{th2}, there are at least
$\maxdeg(\partial \tau)\geq 4$ disjointly eligible tets.
Choose one, say $t$, whose two boundary faces are disjoint from those of
$t_1$.
If removing $t$ leaves one triangulated ball, remove $t$ and continue
ending with $t_1$.
If removing $t$ splits the boundary into two conjoined spheres, then
$\tau-t$ splits into two triangulated balls $\tau_1,\tau_2$ by
Theorem \ref{th1} and Theorem \ref{th2}.
Assume $t_1\in \tau_1$, and let $t_2$ be a tet of $\tau_2$ adjacent to
the gluing face of $t$.
Shuck $\tau_2$ ending with $t_2$, remove $t$, and then shuck
$\tau_1$ ending with $t_1$.
Every removal is the reverse of attaching along one or two faces, so the
reversed order is the desired shelling.
\proofend

\section{Flag complex}

A simplicial complex is a \emph{flag complex}
(or clique complex) if every clique in its $1$-skeleton occurs as a
simplex.
A simplicial triangulation $\sigma$ of $S^2$
is flag exactly when it is prime
(not the connected sum of smaller complexes)
and is not the boundary of a tetrahedron.

We can decompose any $\sigma$ as the connected sum of
prime components $\sigma_i$.
From the point of view of taut fillings these components are independent:
any taut filling of $X(\sigma)$ is the union of taut fillings of the
$X(\sigma_i)$.
We now prove that, regardless of whether $\sigma$ itself is prime,
every taut filling is flag.

\begin{theorem} \label{th4}
If $\sigma$ is a simplicial triangulation of $S^2$ and $M$ is a taut
filling of $X(\sigma)$, then $K(M)$ is a flag complex.
\end{theorem}

\proofstart
Let $\tau = K(M)$.
We must rule out the following minimal failures of flagness:
\begin{enumerate}
\item
the three edges of a triangle $\{A,B,C\}$ without the face $ABC$;
\item
the six edges of the $K_4$ on $\{A,B,C,D\}$ without the tet $ABCD$;
\item
the ten edges of the $K_5$ on $\{A,B,C,D,E\}$.
\end{enumerate}

If Configurations 1 and 2 are absent, so is Configuration 3.  Suppose the
five vertices $\{A,B,C,D,E\}$ carried all ten edges.  With no empty $K_4$,
all five tetrahedra on four of these vertices are simplices of $\tau$; let
$U \subset M$ be the subchain formed by those five tetrahedra.  Each of the
ten triangles on the five vertices is a face of exactly two of the five
tetrahedra, and by the pseudomanifold condition from Theorem \ref{th2} these are the
only tetrahedra of $M$ containing it. Thus the triangle is interior, so it is
not a face of $\sigma$ and has coefficient
$0$ in $\del M = X(\sigma)$.  Hence $\del U$ vanishes on every triangle, so
$\del U = 0$.  By Corollary \ref{nosubcycle} this forces $U = 0$,
contradicting that $U$ contains five tetrahedra.

We rule out Configurations 1 and 2 by minimal counterexample.
If $\sigma$ has a degree-$3$ vertex, then $\sigma$ is a connected sum with
one summand the tetrahedral boundary; removing that summand from $\sigma$ and
the corresponding tetrahedron from $\tau$ gives a smaller counterexample.
So assume $\sigma$ has no degree-$3$ vertex.
By the eligible-tet count of Theorem \ref{th2} there are then at least
$\maxdeg(\sigma)\geq 4$ eligible tetrahedra with pairwise disjoint pairs of
boundary faces; in particular their flip edges are distinct, since the two
boundary faces sharing an edge are the only two faces of $\sigma$ along it.
Removing an eligible tet flips its edge of $\sigma$, leaving either a smaller
triangulated sphere or two spheres conjoined along the new edge.  In the
conjoined case a forbidden $K_3$ or $K_4$, all of whose vertices are joined
by edges, lies on one side, where it gives a smaller counterexample.

\emph{Configuration 1.}  Removing an eligible tet deletes from $\tau$ only
the flipped edge.  The eligible tets have distinct flip edges, and the
forbidden triangle has only three edges, so among at least four eligible
tets some flip edge avoids all three; removing that tet preserves the
triangle, a contradiction.

\emph{Configuration 2.}  Assume Configuration 1 is already ruled out, so the
six edges on $\{A,B,C,D\}$ make all four triangles $ABC,ABD,ACD,BCD$
simplices of $\tau$.  We claim that removing \emph{any} eligible tet $t$
leaves the empty $K_4$ in place, so a single removal already gives a smaller
counterexample.  Removal creates no tetrahedron, so $\{A,B,C,D\}$ is still
not one.  Each of the six edges survives: for an edge, say $AB$, the
triangles $ABC$ and $ABD$ are simplices, and $t$ contains at most one of
them, since containing both would force $\{A,B,C,D\}\subseteq t$, that is
$t=\{A,B,C,D\}$, which is not a tetrahedron.  So one of them, say $ABD$, is
not contained in $t$; being a simplex it is a face of some tetrahedron
$u\neq t$, and $AB\subseteq ABD\subseteq u$ shows $AB$ survives the removal
of $t$.  Hence all six edges remain and the empty $K_4$ persists, giving a
smaller counterexample.  No maximum-degree bound and no octahedral special
case is needed.
\proofend

\bibliography{taut}
\bibliographystyle{amsplain}
\end{document}